\documentclass{amsart}
\usepackage{amssymb}
\usepackage[]{amsthm}
\usepackage[english,polish]{babel} %PS
\usepackage[OT4,plmath]{polski}
\usepackage[cp1250]{inputenc}% Polish letter codding will suppress this lines.

\usepackage{amsmath}
\usepackage{subeqnarray}
\usepackage[colorlinks=true]{hyperref}

\newcommand*{\doi}[1]{\href{http://dx.doi.org/#1}{ doi: #1}}
\newcommand*{\ZBL}[1]{\href{http://www.zentralblatt-math.org/zmath/en/advanced/?q=an:#1&format=complete}{Zbl #1}}

\newtheorem{theorem}{Theorem}
\newtheorem{lemma}[theorem]{Lemma}
\newtheorem{remark}[theorem]{Remark}
\newtheorem{conclusion}{Conclusion}
\newtheorem{conjecture}{Conjecture}
\def\rightbox{\protect\vspace*{-2ex}
\begin{flushright}\(\blacksquare\)\end{flushright}}
\newenvironment{Proof}{{\sc Proof.}\hspace{1mm}}{\rightbox}

\def\fM{{\mathfrak M}}
\def\cF{{\mathcal F}}
\def\cG{{\mathcal G}}
\def\cD{{\mathcal D}}
\def\bP{{\bf P}}
\def\bE{{\bf E}}
\def\one{{\mathbb I}}
\def\bbE{{\mathbb E}}
\def\bS{{\mathbb S}}
\def\vtilde{\tilde v}
\def\Dtilde{\tilde \cD}
\def\stilde{\tilde s}
\def\wtilde{\tilde w}
\newcommand*{\argmax}{\mathop{\mathrm{arg\,max}}\displaylimits}
\title[Duration problem a review]{Duration problem: basic concept and some extensions}
\author[Z. Porosiński]{Zdzisław Porosiński}
\address[Z.~Porosiński]{Wrocław University of Technology \\ \indent Wybrzeże Wyspiańskiego 27 \\ \indent 50-370 Wrocław, Poland}
\email{Zdzislaw.Porosinski@pwr.edu.pl}
\author[M. Skarupski]{Marek Skarupski}
\address[M. Skarupski]{Wrocław University of Technology \\ \indent Wybrzeże Wyspiańskiego 27,\\ \indent 50-370 Wrocław, Poland}
\email{Marek.Skarupski@pwr.edu.pl}
\author[K. Szajowski]{Krzysztof Szajowski}
\address[K. Szajowski]{Wrocław University of Technology \\ \indent Wybrzeże Wyspiańskiego 27 \\ \indent 50-370 Wrocław, Poland}
\email{Krzysztof.Szajowski@pwr.edu.pl}

\begin{document}
%\vspace{-5.5ex}
%\Poczatek
%\Chapter
%\pagenumbering{roman}
\setcounter{page}{1} %%This command starts the numerations of pages
\selectlanguage{english}%\Polskifalse
%\selectlanguage{polish}\Polskitrue
%\setcounter{page}{1} %%This command starts the numerations of pages

\begin{abstract}
We consider a sequence of independent random variables with the known distribution observed sequentially. The observation $n$ is assumed to be a value of one order statistics such as $s:n$-th, where $1 \leq s\leq n$. It the instances following the $n$th observation it may remain of the $s:m$ or it will be the value of the order statistics $r:m$ (of $m> n$ observations). Changing the rank of the observation, along with expanding a set of observations there is a random phenomenon that is difficult to predict. From practical reasons it is of great interest. Among others, we pose the question of the moment in which the observation appears and whose rank will not change significantly until the end of sampling of a certain size. We also attempt to answer which observation should be kept to have the ``good quality observation'' as long as possible. This last question was analysed by Ferguson, Hardwick and Tamaki (1991) in the abstract form which they called the \emph{problem of duration}.

This article gives a systematical presentation of the known duration models and a new modifications. We collect results from different papers on the duration of the extremal observation in the no-information (denoted as rank based) case and the full-information case. In the case of non-extremal observation duration models the most appealing are various settings related to the two extremal order statistics. In the no-information case it will be the maximizing duration of owning the relatively best or the second best object. The idea was formulated and the problem was solved by Szajowski and Tamaki (2006). The full-information duration problem with special requirement  was presented by Kurushima and Ano (2010).
\nocite{KurAno10:FIDP,SzaTam:arXiv0902.0232S}
\end{abstract}
\maketitle
%*******************************************************************
%DO NOT FORGET TO RESET THE EQUATION COUNTER TO 0 AT THE HEAD OF EACH SECTION
%*******************************************************************
%*******************************************************************
%FIRST SECTION - DO NOT FORGET THE SETCOUNTER{EQUATION}{0}
%*******************************************************************
%\vspace{-5.5ex}
\setcounter{equation}{0}
\section{Introduction}%\input{sec0.tex}
%---Introduction---%
It was Ferguson, Hardwick and Tamaki~\cite{ferhartam92:own} who formulated the duration problem as the optimal prediction the relative extremal observation keeping the leading position for the highest period. The basic formulation was for the classical no-information secretary problem. It is a sequential selection problem which is a variation of the classical secretary problem (CSP) treated, for example, by Gilbert and Mosteller~\cite{gilmos66}. The aim of the CSP is to examine items ranked from $1$ to $N$ by the random selection without replacement, one at a time, and to win wchich means to stop at any item whose overall (absolute) rank belongs to the given set of ranks (in the basic version this set contains the rank $1$ only), given only the relative ranks of the items drawn so far. Since the articles by Gardner~\cite{gar:sci60} the secretary problem has been extended and generalized in many different directions. Excellent reviews of the development of this colourful problem and its extensions have been given by Rose~\cite{rose}, Freeman~\cite{fre} and Samuels~\cite{sam91:secretary}. The deepest analysis of the assumptions and their consequences was made by Ferguson~\cite{ferhartam92:own}.

\subsection{\label{intro_1}Duration problems for the no-information case.}%\input{sec0_1.tex}
The basic form of the duration problem can be described as follows. A set of $N$ rankable objects appears as in CSP. As each object appears, we must decide whether to select or reject it on the basis of the relative ranks of the objects. The payoff is the length of time we are in possession regarding a relatively best object. Thus we will only select the relatively best object, receiving a payoff of one as we do so and an additional one for each new observation as long as the selected object remains relatively best.

Ferguson, Hardwick and Tamaki~\cite{ferhartam92:own} considered various duration models in quite some detail. Moreover, they mention that they had discussed the duration problem also for a random number of arrivals in continuous time. The solution is indeed easy if the arrival process is the Pascal process, since then (see Bruss and Rogers~\cite[Theorem 2]{brurog91:pascal}), the corresponding record arrival process is Poisson. The latter gives then also access to Poisson embedding (Bruss and Rogers~\cite{brurog91:Poisso}) displaying an interesting duality. However, Ferguson et al.~\cite{ferhartam92:own}  confined themselves throughout their study to the duration problem for the relatively best items. 
%The list of no-information case problems considered there is presented in the section~\ref{sec2shell}.
In this paper, we attempt to extend the problems to choose and keep the items having the relatively high leading position for the long period. As a simple example we refer to a relatively best or the second best object as a candidate. The focus is the case where each time we receive a unit payoff as long as  either of the chosen objects remains a {\bf candidate}. Obviously only candidates can be chosen, the objective being to maximize the expected payoff. This problem can be viewed from another perspective as follows. Let us observe at moment $i$ the relatively second  candidate and let us denote $T(i)$ the time  of the first candidate after time $i$ (\emph{i.e.} the relatively best or the second best item) if there is one, and $N+1$ if there is none. If we observe at $i$ the relatively best item then $T(i)$ is the moment when new item appears, which changes the relative rank of $i$th item to the non-candidate rank. The time $T(i)-i$ is called duration of the candidate selected at time $i$. The objective is to find a stopping time $\tau^{*}$ such that
\begin{equation}\label{problem}
 v_N=E\left[ \frac{T(\tau^*)-\tau^*}{N}\right]
 =\sup\limits_{\tau\in \fM^N} E\left[\frac{T(\tau)-\tau}{N}\right],
\end{equation}
where $\fM^N$ denotes the set of all stopping times.
\par
\subsection{\label{intro_2}Duration problems for the full-information case.}
The second group of models in \cite{ferhartam92:own} are those related to the full-information best choice problem. We observe sequentially i.i.d. random variables from the known distribution. Without the loss of generality we can assume that they come from the uniform distribution on the interval $[0,1]$. Suppose that we want to maximize the time in which the selected object maintains its quality (e.g. being the relatively best one). In many cases we see that the reward is related with the win probability. A typical optimal stopping problem in full information case was first studied by Gilbert and Mosteller \cite{gilmos66}. The comparison of results was derived by Gnedin \cite{Gne05:Objectives}.\\

The paper is organized as follows.  In Section~\ref{sec1} the solution of the problems formulated in the section~\ref{intro_1} is presented. A Markov chain optimal stopping problem equivalent to the duration problem and the optimal strategy will be formulated in Section~\ref{sec1.1} and derived in Section~\ref{sec1.2}. This section is based mainly
on the suggestion from~\cite{dynyush} and the results by~\cite{sza82:ath} and \cite{sucsza02:no-information}. In Section~\ref{sec1.3} the problem of stopping on the relatively best or second best is shown. The discounted no-information duration problem is described in in Section~\ref{sec1.4}. In the Section~\ref{sec2} we consider full information duration problem (FIDP) in finite horizon with and without recall possibility. The duration problem is transformed to the optimal stopping problem for the Markov process similar to this applied by Bojdecki~\cite{boj79:probability} and Prosi\'nski and Szajowski~\cite{porsza90:selection} (see Section~\ref{sec2.1}). In Section \ref{sec2.2} we present results related to the best-choice duration problem (BCDP) with and without recall possibility. In section \ref{sec2.3} we switch into random horizon full information duration problem (RHFIDP). In Section \ref{sec2.4} there is a description of the full information duration problem with the unbounded horizon. We extend FIDP into duration of owning the relatively best or second best object in Section \ref{sec2.5}. We show some results for the unbounded horizon. Various directions of the extensions based on the main idea of the duration problem are the subject of the conclusion section. 
%\vspace*{0.5cm}
%*******************************************************************
%FIRST SECTION - DO NOT FORGET THE SETCOUNTER{EQUATION}{0}
%*******************************************************************
%\vspace*{0.6cm}

\section{\label{sec1}Markov model for no-information duration problems}
The models which are considered in this study are the so called no information models where the  decision to select an object is based only on the relative ranks of the objects observed so far. Let $\bS=\{1,2,\ldots,N\}$ be the set of ranks of items $\{x_1,x_2,\ldots,x_N\}$ and  $\{X_1,X_2,\ldots,X_N\}$  their permutation. We observe sequentially the permutation of items from the set $\bS$. The mathematical model of such an experiment is the probability space $(\Omega,\cF,\bP)$. The elementary events are permutations of the elements from $\bS$ and the probability measure $\bP$ is the uniform distribution on $\Omega$. The observation of random variables $Y_k$, $k=1,2,\ldots,N$, generate the sequence of $\sigma$-fields $\cF_k=\sigma\{Y_1,Y_2,\ldots,Y_k\}$, $k=1,2,\ldots,N$. The random variables $Y_k$ are independent and $\bP\{Y_k=i\}=\frac{1}{k}$. 

We assume that all permutations are equally likely. If $X_k$ is the rank of $k$-th candidate we define
\begin{equation}\label{relrank1}
Y_{kj}=\#\{1\leq i\leq j: X_i\leq X_k\}%Y_k=\#\{1\leq i\leq k: X_i\leq X_k\}.
\end{equation}
the running rank  of $k$-th object at the moment $j\geq k$. The random variable $Y_k=Y_{kk}$ is called {\it relative rank} of $k$-th candidate with respect of the items investigated to the moment $k$. Let $A\subset \bS$. Subsequently, the appearance of the next candidate after $i$ is the moment $S_A(i)=\inf\{k>i:Y_k\in A\}$ and the maturity $T_r(i)=\inf\{k\geq i:\bP\{X_k\in A|Y_{ik}=r\}=0\}$ of the candidate with the relative rank $r$ at moment $i$ is defined. The maturity can be represented in the equivalent manner as $T_r(i)=\inf\{k\geq i: Y_i=r, \ Y_k\in A,\ Y_{ik}\notin A\}$.
\begin{remark} 
If $A=\{1,2,\ldots,s\}$ and
\begin{subequations}\label{Durat}
\begin{align}
\label{durat1}
S_A(i)&=\inf\{k>i:Y_k\in A\}\\
\nonumber \intertext{then $T_s(i)=S_A(i)$, and for any $r\in A$}
%\nonumber\mbox{and for any $r\in A$\hspace{5em}}&&\\ %\hspace{11em}
\label{durat3}
T_r(i)&=\inf\{k\geq i:Y_{ik}\not\in A\}.
\end{align}
\end{subequations}
For $r<s$ the duration of the candidate with rank $r$ at the moment $i$ is dependent on the items appearing between $i$ and $k$. Changing the rank of a candidate when a new candidate approaches, does not always mean that one ceased to be a candidate. 
\end{remark} 

Denote by $\fM^N$ the set of all Markov moments $\tau$ with respect to $\sigma$-fields $\{\cF_k\}_{k=1}^N$. The decision maker observes the stream of relative ranks. When $Y_i\in A=\{1,2,\ldots,s\}$ it is the potential candidate for the absolutely $r$th item, $r\in A$. Sometimes it is enough to keep such a candidate for a period of time to get profit which is proportional to the shell file of a candidate (the second kind of duration of the candidate). The random variable $T(i)$ is defined as the moment when the keeping candidate stops to be the candidate (the maturity of the candidate). Let us consider the possibility of recall force to think about the rejected candidates. We define $\delta_r(i)=\sup_{1\leq j\leq i}\{Y_j=r\}$ \emph{the actual position of the relative $r$ at the moment $i$}. $\delta_r(i)$ is the random variable measurable with respect to $\cF_i$. The recall option means the possibility of returning to the last candidate who has the relative rank $r^\star=\mbox{arg}\max_{\{s\in A\}}\delta_s(i)$ or to the candidate with the given rank, \mbox{e.g} $r=1$. The present history at $i$ and $Y_{k,s}$ for $k\in\{i,\ldots,N\}$, $s=i,i+1,\ldots,N$ allows to define the maturity $T(i)$ for various models. In the next part examples of various definitions of maturity will be shown, and therefore different definitions of duration will be presented.
  
\subsection{\label{sec1.1}The optimal stopping problem for the embedded Markov chain}
  %\input{sec1.1.tex}
%%%%%%%%%%%%%%%%%
\subsubsection{Embedded Markov chain.} Let $a=\max(A)$. The function $\varphi(k,r)$ defined in \eqref{varphi1} is equal to
$0$ for $r>a$ and non-negative for $r\leq a$. It means that it is rational to choose an item for keeping
at the moment $k$ when the state $(k,r)$ such that $r\leq a$.
Define $W_0=(1,Y_1)=(1,1)$, $\gamma_t=\inf\{r>\gamma_{t-1}:Y_r\leq
\min(a,r)\}$ ($\inf\emptyset=\infty$) and $W_t=(\gamma_t,Y_{\gamma_t})$.
If $\gamma_t=\infty$ then define $W_t=(\infty,\infty)$. $W_t$ is  the
Markov chain with the state space $\bbE=\{(s,r):s\in\{1,2,\ldots,N\},r\in A\}\cup\{(\infty,\infty)\}$
and the following transition probabilities
(see~\cite{sza82:ath})
\begin{equation}\label{trprob}
  \begin{split}
  p(r,s)&=\bP\{W_{t+1}=(s,l_s)|W_t=(r,l_r)\} \\
    & =\begin{cases}
      \frac{1}{s},             &\text{ if $r<a$, $s=r+1$},\\
      \frac{(r)_a}{(s)_{a+1}}, &\text{ if $a\leq r< s$},\\
      0,                       &\text{ if $r\geq s$ or $r<a$, $s\neq r+1$,}
    \end{cases}
  \end{split}
\end{equation}
with  $p(\infty,\infty)=1$, $p(r,\infty)=1-a\sum_{s=r+1}^Np(r,s)$, where
$(s)_a=s(s-1)(s-2)\ldots(s-a+1)$, $(s)_0=1$. We denote $T\varphi(k,r)=\bE_{(k,r)} \varphi(W_1)$
the mean operator for the function $g:\bbE\rightarrow \Re$. Let
$\cG_t=\sigma\{W_1,W_2,\ldots,W_t\}$ and $\tilde{\fM}^N$ be the set of stopping
times with respect to $\{\cG_t\}_{t=1}^N$. Since $\gamma_t$ is increasing, then
we can define $\tilde{\fM}^N_{r+1}=\{\sigma \in \tilde{\fM}^N : \gamma_\sigma>r\}$.

Let $\bP_{(k,r)}(\cdot)$ be the probability measure related to the Markov chain
$W_t$, with the trajectory starting in the state $(k,r)$ and $\bE_{(k,r)}(\cdot)$ the
expected value with respect to $\bP_{(k,r)}(\cdot)$. From \eqref{trprob} we can
see that the transition probabilities do not depend on relative ranks, but only
on the moments $k$ where the items with the relative rank $r\leq \min(a,k)$ appear. Based
on the following lemma we can solve the problem~\eqref{problem} with the gain
function \eqref{varphi} using the embedded Markov chain $(W_t,\cG_t,\bP_{(1,1)})_{t=0}^N$.
\begin{lemma} (see \cite{sza82:ath})
  \begin{equation}\label{basic1}
  \bE w_N(k+1,Y_{k+1})=\bE_{(k,r)} w_N(W_1) \mbox{ for every $r\leq \min(a,k)$.}
  \end{equation}
\end{lemma}

%%%%%%%%%%%%%%%%%%%%%%%%%%%%%%%%%%%%%%%%%%%%%%%%%%%%%%%%%%%%%%%%%%%%%%%%%%
\subsubsection{\label{sec22s}The optimal stopping problem.} Let $T(i)=\zeta_iT_1(i)$ and $\xi_i$ is an additional restriction (the requirement on the chosen item). The aim is to find  $\tau^\star\in\fM^N$ such that: 
\begin{equation}\label{FHD_BCPwr}
\bE\left[\frac{T(\tau^\star)-\tau^\star}{N}\xi_{\tau^\star}\right]=\sup_{\tau\fM^N}\bE\left[\frac{T(\tau)-\tau}{N}\xi_\tau\right].
\end{equation}

Let us observe that for any $\tau\in\fM^N$
\begin{eqnarray*}\label{problem1}
 \bE\left[\frac{T(\tau)-\tau}{N}\xi_\tau\right]
% &=& \sum_{r=1}^N\int_{\{\tau=r\}}\frac{T_r-r}{N}\xi_r d\bP=
%    \sum_{r=1}^N\int_{\{\tau=r\}}\bE\{\frac{T_r-r}{N}\xi_r|\cF_r\}d\bP\\
 &=& \sum_{r=1}^N\int_{\{\tau=r\}}\bE\{\frac{T_r-r}{N}\xi_r|Y_r\}d\bP
  = \bE \varphi(\tau,Y_\tau).
\end{eqnarray*}

%%%%%%%%%%%%%%%%%%%%%%%%%%%%%%%%%%%%%%%%%%%%%%%%%%%%%%%%%%%%%%%%%%%%%%%%%%
\subsubsection{Recursive algorithm} Let $\fM^N_r=\{\tau\in\fM^N: r\leq \tau\leq N\}$ and
$\wtilde_N(r)=\sup_{\tau\in\fM^N_r}\bE \varphi(\tau,Y_\tau)$. The following algorithm
allows to construct the value of the problem $v_N=w_N(1,1)$.
\begin{equation}\label{step1}
  \wtilde_N(N)=\bE \varphi(N,Y_N)%=\frac{2}{N}.
\end{equation}
Let
\begin{subeqnarray}\label{steps}
\slabel{step2}
  w_N(N,r)&=&\left\{
            \begin{array}{ll}
                  1, &\mbox{ if $r\in A$},\\
                  0, &\mbox{ otherwise,}
            \end{array}
                  \right.\\
\slabel{step3}%
  w_N(k,r)&=&\max\{\varphi(k,r),\bE w_N(r+1,Y_{r+1})\},\\
\slabel{step4}%
  \wtilde_N(k)&=&\bE w_N(k,Y_k)=\frac{1}{k}\sum_{r=1}^k w_N(k,r).
\end{subeqnarray}
We have then $v_N=\wtilde_N(1)$. The optimal stopping time $\tau^*$ is defined
as follows: one has to stop at the first moment $k$ when $Y_k=r$, unless
$w_N(k,r)>\varphi(k,r)$. We can define the stopping set as $\Gamma=\{(k,r): \varphi(k,r)\geq \wtilde_N(k)\}$.
  
\subsection{\label{sec1.2}Classical no-information BC duration problem}
It is not difficult to formalize the duration problem for BC as with  recall as without recall and also when the additional requirement concerning the absolute rank of the selected object is added. In order to present the problem of the duration time for BCP considered in \cite{ferhartam92:own} we assume that $A=\{1\}$, $\zeta_n(\omega)=\one_{\{Y_n\in A\}}(\omega)$ and $\zeta_n^\star(\omega)=\one_{\{X_n\in A\}}(\omega)$. 

\subsubsection{Finite horizon duration problem of BCP without recall (\cite[Sec. 2.2]{ferhartam92:own}).} Let $T(i)=\zeta_iT_1(i)$ and $\xi_i=1$. The aim is to find  $\tau^\star\in\fM^N$ such that: 
\begin{equation}\label{FHD_BCPwor}
\bE\left[\frac{T(\tau^\star)-\tau^\star)}{N}\xi_{\tau^\star}\right]=\sup_{\tau\fM^N}\bE\left[\frac{(T(\tau)-\tau)}{N}\xi_\tau\right].
\end{equation}
It is the first setting of the problem. In ~\cite{ferhartam92:own} the authors observed that the pay-offs in the problem for the threshold rules are exactly there same as the pay-offs for the threshold rule for the best choice secretary problem with an unknown, random number of options that has the uniform distribution on $\bS$ (see the results by Presman and Sonin~\cite{preson72:random}, Rasmussen and Robbins\cite{rasrob75} and a general method of Samuels~\cite{sam91:secretary} showing the relation of the random horizon problems to the problems with cost ). The single threshold $r^\star_N$strategy is optimal having the asymptotic $\lim_{N\rightarrow\infty}\frac{r^\star_N}{N}\cong \mbox{e}^{-2}$ and the problem value $2\mbox{e}^{-2}$. 
  
\subsubsection{Finite horizon duration problem of BCP with recall (\cite[Sec. 2.2]{ferhartam92:own}).} Let $T(i)=T_1(i)$ and $\xi_i=1$. The aim is to find  $\tau^\star\in\fM^N$ in the problem of \eqref{FHD_BCPwr} with this new definition of the maturity moment. This second setting of the problem has the solution which has a simple relation with the solution of BCP. Namely, if $k_N^\star$ is the optimal threshold for BCP then the optimal threshold for the duration problem of BCP with recall is $K_N=r_N^\star-1$ for $N\geq 2$. It is also the optimal rule for the BCP with an unknown, random number of options having the uniform distribution on $\bS$ and possible recall.  

\subsubsection{Duration problem without recall with choice of the best (\cite[Sec. 2.3]{ferhartam92:own}).} Let us define the maturity moment as $T(i)=\zeta_i T_{1}(i)$ and $\xi_i=\zeta_i^\star$. The positive pay-off is when $T(i)=N+1$ only. The expected fraction of duration $v_k$ for the threshold strategy $k$ is equal $v_k=\frac{N+1-k}{N}\bP\{T(k)=N+1\}=\frac{N+1-k)k}{N^2}$ which is unimodal with mode at $L_N=\left\lfloor \frac{N+1}{2}\right\rfloor$. Thus the optimal rule is among the threshold rules with $k$ lower than $L_N$. The asymptotic $\alpha=\lim_{N\rightarrow\infty}\frac{r^\star_N}{N}\cong 0.20388$ is the solution of the equation $-\log(x)-2+2x=0$. The limiting value of the expected pay-off is $0.1618$.

\subsubsection{Duration problem with recall and choice of the best (\cite[Sec. 2.3]{ferhartam92:own}).} Let us define the maturity moment as $T(i)=T_1(i)$ and $\xi_i=\zeta_i^\star$. The positive pay-off is when $T(i)=N+1$ only. The optimal rule is the fixed sample size rule that stops at $L_N$. The asymptotic optimal return is $0.25$.
  
\subsection{\label{sec1.3}The duration of the best or the second best}
\subsubsection{\label{TiDistribution}Distribution of a maturity moment.} Let $A=\{1,2\}$.    The model without recall is considered. The maturity of the candidate at $i$ is equal: 
\begin{equation}\label{maturity1or2}
T(i)=\sum_{r=1}^2\one_{\{Y_{\delta_r(i)}=r\}} T_{r}(i),
\end{equation}
where $\delta_r(i)$, the position of the relatively $r$th at the moment $i$, is equal $i$ for $r=1,2$. The conditional distribution of $T(i)$ is following:
\begin{description}
\item{$\mathbf{Y_i=2:}$} In this case $T(i)=k$ when $Y_i=2,Y_{i+1}>2,Y_{i+2}>2,\ldots,Y_{k-1}>2,Y_k\in A$.
We have for $i<k\leq N$:

%\begin{subequations}\label{rel2timed}
%\bP\{T(i)=k|Y_i=2\}=\left\{\begin{align}\slabel{rel2timedk}
%&=\frac{2(i-1)i}{(k-2)(k-1)k}&\text{for $ i<k \leq N$;}\\
%\slabel{rel2timedN}
%&=1-\sum_{s=i+1}^N\frac{2(i-1)i}{(s-2)(s-1)s} =\frac{i(i-1)}{N(N-1)}&\text{ for $k=N+1$.}
%\end{align}\right.
%\end{subequations}
\begin{subequations}\label{rel2timed}
\begin{align}\slabel{rel2timedk}
\bP\{T(i)=k|Y_i=2\}&=\frac{2(i-1)i}{(k-2)(k-1)k};\\
\slabel{rel2timedN}
\begin{split}
\bP\{T(i)=N+1|Y_i=2\}&=1-\sum_{s=i+1}^N\frac{2(i-1)i}{(s-2)(s-1)s} \\
&=\frac{i(i-1)}{N(N-1)}.
\end{split}
\end{align}\end{subequations}

%\begin{subeqnarray}\label{rel2timed}
%\begin{align}\slabel{rel2timedk}
%\bP\{T(i)=k|Y_i=2\}&=&\frac{2(i-1)i}{(k-2)(k-1)k};\\
%\slabel{rel2timedN}
%\bP\{T(i)=N+1|Y_i=2\}&=&1-\sum_{s=i+1}^N\frac{2(i-1)i}{(s-2)(s-1)s} =\frac{i(i-1)}{N(N-1)}.
%\end{align}\end{subeqnarray}
%\begin{equation}
% u(x) =
%  \begin{array}{ll}
%   \exp{x} & \text{if $ x \geq 0$} \\
%   1       & \text{if  $x < 0$}
%  \end{array}
%\end{equation}

%%%%%%%%%%%%%%Rel 1 %%%%%%%%%%%%%%%%%%%%%%%%%%%%%%
\item{$\mathbf{Y_i=1; T(i)\leq N:}$} the random variable $T(i)=k$ if there exists $s\in\{i+1,\ldots,k-1\}$ such that
 $Y_i=1,Y_{i+1}>1,Y_{i+2}>1,\ldots,Y_{s-1}>1,Y_s=1,Y_{s+1}>2,\ldots,Y_{k-1}>2,Y_k\in A$. We have for $i<k\leq N$
\begin{subeqnarray}\label{rel1timed}
\slabel{rel1timedk}
\bP\{T(i)=k|Y_i=1\}&=&\frac{2i(s-i+1)}{(k-2)(k-1)k}.
\end{subeqnarray}

\item{$\mathbf{Y_i=1,T(i)=N+1,Y_{kN}=1 \textrm{ or } 2:}$} We have 
\begin{subeqnarray}\label{rel1timed1}
\slabel{rel1timedN1}
\bP\{T(i)=N+1,Y_{iN}=1|Y_i=1\}&=&1-\sum_{s=i+1}^{N}\frac{i}{(s-1)s}\\
\nonumber&&=\frac{i}{N};\\
\slabel{rel1timedN2}
\bP\{T(i)=N+1,Y_{iN}=2|Y_i=1\}&=&
\sum_{s=i+1}^{N}\frac{i}{(s-1)s}\\
\nonumber&&(1-\sum_{k=s+1}^N\frac{2(s-1)s}{(k-2)(k-1)k})\\
\nonumber&&=\frac{i(N-i)}{N(N-1)}\\
\slabel{rel1timedN1i2}
\bP\{T(i)=N+1|Y_i=1\}&=&\frac{i}{N}+\frac{i(N-i)}{N(N-1)}\\
\nonumber&&=\frac{i(2N-i-1)}{N(N-1)}
\end{subeqnarray}
\end{description}

\begin{remark}
The solution of the problem \eqref{problem} with $T(i)$ given by \eqref{maturity1or2} will be performed by its change to the optimal stopping problem for the embedded Markov chain. In the case without recall there are no additional restrictions and $\xi_n=1$ for $n=1,2\ldots,N$. However, there is the obvious and interesting problem of taking into account the value of a candidate who is kept until moment $n$ (see \cite{Gne07:Rank_depend}). Especially, when the value of the candidate is changed over time. 
\end{remark}

\begin{remark}
 The possibility of recall requires an additional clarification. The natural models are as follows:
\begin{description}
\item[(i)] the possibility of returning to the best so far candidate means $T(i)=T_{1}(i)$;
\item[(ii)] the possibility to return to the last candidate is defined as $T(i)=\one_{\{Y_{\delta_{r^\star}(i)}=r^\star\}} T_{r^\star}(i)$, where $r^\star=\argmax_{\{s\in A\}}\delta_s(i)$.
\end{description}  
Both approaches require access to the history of the observed random variables during the selection process, each at a different extent. 

The duration problem which requires the selected object to be of the prescribed absolute rank (which is not higher than $a$) and matching the expectations also needs further investigation. 

%The duration problem with requiring that the selected object be of the prescribed absolute rank also requires additional clarification. These models will be subject of further research.   
\end{remark}

%%%%%%%%%%%%%%%%%%%%%%%%%%%%%%%%%%%%%%%%%%%%%%%%%%%%%%%%%%%%%%%%%%%%%%%%%%
\subsubsection{\label{sec22}The optimal stopping problem for the best or the second best duration problem.} In the following lemma the function $\varphi(\cdot)$ is calculated. The final form of it is using 
the digamma function ($\digamma$-function) $\psi_n(z)$ (see Abramowitz and Stegun \cite{abrste72:functions} p. 260).
For $n=0$ we will use the denotation $\psi(z)$.
This function is defined as $n$th logarithmic derivative of the Euler gamma function $\Gamma(z)$
\begin{eqnarray*}
\psi_n(z)&=&\frac{d^{n+1}}{d\!z^{n+1}}\ln\Gamma(z)= \frac{d^n}{d\!z^n}\frac{\Gamma^{'}(z)}{\Gamma(z)}.
\end{eqnarray*}

\begin{lemma}\label{thepayoff}
The payoff function $\varphi(k,r)$ has the form
{\renewcommand{\arraystretch}{2}
\begin{equation}\label{varphi}
\varphi(k,r)=\left\{\begin{array}{ll}
\frac{k}{N^2}\left(1+k-N-2N\psi(k)+2N\psi(N)\right)&\mbox{ for $r=1$,}\\
\frac{k(N-k+1)}{N^2}&\mbox{ for $r=2$,}\\
0&\mbox{ otherwise.}
\end{array}
\right.
\end{equation}
}
\end{lemma}
\begin{Proof}
Based on the distribution of the random variable $T(k)$ and the equality $\psi(p+1)-\psi(p)=\frac{1}{p}$ for
the digamma function  we get

\begin{subeqnarray}
\label{varphi1}
\slabel{varphi1A}
\varphi(k,1)&=&\bE\{\frac{T(k)-k}{N}|Y_k=1\}\\
\nonumber&&=\frac{k}{N^2}\left(1+k-N-2N(\psi(k)-\psi(N))\right)\\
\slabel{varphi1B}
\varphi(k,2)&=&\bE\{\frac{T(k)-k}{N}|Y_k=2\}=\frac{k(N-k+1)}{N^2}
\end{subeqnarray}
\end{Proof}

\subsubsection{Solution of the optimal shelf life problem for the best and the second best.} First of all the form of $T\varphi(k,r)$ for $(k,r)\in \bbE$ will be given.
\begin{lemma}\label{exppayoff}
The expected payoff of the function $\varphi(\cdot)$ with relation to the embedded Markov chain
$\{W_t\}$ has the following form:
\begin{equation}\label{basic2}
  T\varphi(k,r)=\frac{(N-k)((2N-1)k+N-1)}{N^2(N-1)}+2\frac{k}{N^2}(\psi(N)-\psi(k)).
  \end{equation}
\end{lemma}
\begin{Proof}
The definition of the embedded Markov chain \eqref{trprob} and the payoff function $\varphi(\cdot)$ in
the lemma \ref{thepayoff} give
\begin{eqnarray*}
T\varphi(k,r)&=&\sum_{j=k+1}^N\sum_{r=1}^{2}p(k,j)\varphi(j,r)\\
             &=&\sum_{j=k+1}^N\frac{k(k-1)}{j(j-1)(j-2)}\Big(\frac{j(2N(\psi(N)-\psi(j))+N-j-1)}{N^2}\\
             \nonumber&&+\frac{j(N-j+1)}{N^2}  \Big)  \\
             &=&\sum_{j=k+1}^N\frac{k(k-1)}{j(j-1)(j-2)}
             \left(\frac{j(N-j-1)}{N^2}(j+N+\frac{2N}{j}-\frac{2N}{N-1})  \right)\\
             &=&\frac{(N-k)((2N-1)k+N-1)}{N^2(N-1)}+2\frac{k}{N^2}(\psi(N)-\psi(k)).
\end{eqnarray*}
\end{Proof}

Let us denote $A_k(r)=\{(s,r):s>k\}$.
\begin{theorem}\label{expectedpayoff}
There are constants $k_1^\star$ and $k_2^\star$ such that the optimal stopping time for the problem \eqref{problem}
has the form
\begin{equation*}
      \tau^*=\inf\{t:W_t\in A_{k_1^\star}\cup A_{k_2^\star}\}.
\end{equation*}
The value function is 
\begin{eqnarray*}
\vtilde_N(k_1^\star,k_2^\star)&=&
\frac{(N (3N-4)-3)+k_1^\star(N-3)\psi(k_1^\star)}{(N-1)N}\\
			&&\mbox{} \frac{ +k_1^\star(2 (N^2-1)(\psi_{1}(k_2^\star+1)}{(N-1) N}-\frac{ \psi_{1}(k_1^\star + 1)))}{(N-1) N} \\
           &&\mbox{} +\frac{k_1^\star \left(2 (N-1) \psi(N)+(5-3 N) \psi(k_2^\star)\right)}{(N-1) N }\\
&&\mbox{}-\frac{k_1^\star\left(3 N^3+(2 k_2^\star-3) N^2-2 %
          \left({k_2^\star}^2+k_2^\star+2\right)N+{k_2^\star}^2+k_2^\star\right)}%
   {(N-1) N^2 k_2^\star}
\end{eqnarray*}
\end{theorem}
\begin{Proof} The payoff function $\varphi(\cdot,r)$ for $r\in A$ is unimodal. It can be seen following the analysis of the differences $\varphi(k+1,1)-\varphi(k,1)$ which are decreasing when $k\leq N-1$. The comparison  of events related to $T(k)=j$ on $Y_k=1$ and $Y_k=2$ leads to the conclusion that $\varphi(k,1)\geq \varphi(k,2)$ for $k\in\{1,2,\ldots,N\}$. The value function $\wtilde_(k)$ is no-increasing due to the fact of decreasing number of stopping times in $\fM^N_k$. At $k=N-1$ both payoff functions are greater than $\wtilde_N(N-1)$. Let us assume $k_2^\star=\inf\{1\leq k\leq N: T\varphi(k,i)\leq \varphi(k,2)\}-1$. We have for $k>k_2^\star$ and $r=1,2$ that $w_N(k,r)=\varphi(k,r)$ and $\wtilde_N(k)=T\varphi(k,r)$. Let us denote $k_1^\star=\inf\{1\leq k\leq k_2^\star: \wtilde_N(r)<\varphi(k,1) \}$, where $\wtilde_N(k)=\vtilde(k,k_2^\star)$ and for $k<s$ we have
\begin{eqnarray*}
\vtilde_N(k,s)&=&\sum_{j=k+1}^s\frac{k}{j(j-1)}\varphi(j,1)+\frac{k}{s}\wtilde_N(s)\\
  &=&\frac{(N (3 N-4)-3)+k (N-3) \psi(k)   } {(N-1) N}\\
  &&\mbox{} +\frac{k \left(2 (N-1) \psi(N)+(5-3 N) \psi(s)\right)}{(N-1) N}\\
   &&\mbox{}+\frac{ k\left(2 (N^2-1)\left(\psi_1(s+1)-\psi_1(k+1)\right)\right) }{(N-1) N}\\
   &&\mbox{}-\frac{k \left(3 N^3+(2 s-3) N^2-2 \left(s^2+s+2\right) N+s^2+s\right)}%
   {(N-1) N^2 s}.
\end{eqnarray*}
\end{Proof}
  
\subsection{\label{sec1.4}Discounted no-information duration problem (DNIDP)}
	%\input{sec1.4.tex}
%sec2shellvMS1b.tex
The discounted model for no-information duration problem was formulated by Ferguson et a.~\cite{ferhartam92:own}. Methods of discounting applied in their paper assume that the horizon is infinite but the future payoffs are discounted by coefficient $\beta\in (0,1)$ in such a way that an amount of $1$ received at the moment $k$ is worth $(1-\beta)\beta^{k-1}$ at moment $0$. If the $k$-th object is relatively the best, and if the decision maker selects it, the payoff is equal to the discounted maturity moment of the object chosen. The conditional expected return is then
\[
\bE(1-\beta)\beta^{T_1(k)-1} T_1(k)=\sum_{j=k+1}^\infty(1-\beta)\beta^{j-1}j p(k,j)=(1-\beta)k\sum_{j=k}^\infty\frac{\beta^j}{j}.
\] 
Hence $\bE(1-\beta)\beta^{T_1(k)-1} T_1(k)=(1-\beta)\sum_{j=k}^\infty\beta^j\bP(T_1(k)>j)$ because $\bP(T_1(k)>j)=\sum_{s=j+1}^\infty\frac{k}{s(s-1)}=\frac{k}{j}$. 
The problem is to find $\tau^\star$ such that 
\begin{equation}\label{DNIDP}
\bE(1-\beta)\beta^{T_1(\tau^\star)-1} T_1(\tau^\star)=\sup_{\{\tau\in\mathfrak{S}\}}\bE(1-\beta)\beta^{T_1(\tau)-1} T_1(\tau)
\end{equation}
The unimodality of the conditional expected return allows to formulate the optimal stopping time maximizing the expected discounted duration in model formulated in this way. 
\begin{theorem} (Ferguson et al. \cite{ferhartam92:own} ) 
In the discounted duration problem \eqref{DNIDP}, the optimal rule has the threshold form defined by 
\begin{equation}
r^\star(\beta)=\min\{r\geq 1:\sum_{j=r+1}^\infty\frac{\beta^j}{j}\sum_{k=r+1}^j\frac{1}{k-1}\leq \sum_{j=r}^\infty\frac{\beta^j}{j}\}.
\end{equation}
\end{theorem}

The solution of the problem is closely related to the BCP with the random, geometric number of objects. 
	
%*******************************************************************
%SECOND SECTION - DO NOT FORGET THE SETCOUNTER{EQUATION}{0}
%*******************************************************************
%\vspace*{-1ex}

\section{\label{sec2}Markov model for full information duration problem.}
In full information duration problem our aim is to chose the object and hold it as long as it is in its position, and we may base our choice of the stopping time on the true values of the object.

We introduce a Markovian approach (cf. Bojdecki~\cite{boj79:probability}). Let $N \in \mathbb{N}$ be a fixed number and $\lbrace X_{n} \rbrace_{n=1}^{N}$ be a sequence of i.i.d random variables uniformly distributed on the interval $[0,1]$. For $n=1,...,N$ define $\mathcal{F}_{n} = \sigma (X_{1},...,X_{n})$ and let $\mathcal{T}$ denote a set of all stopping times with respect to the family $\lbrace \mathcal{F}_{n} \rbrace_{n=1}^{N}$. Let $\mathcal{T}_{0}$ denote a set of all stopping times $\tau \in \mathcal{T}$ such that
\begin{displaymath}
X_{n} = \max \lbrace X_{1},...,X_{n} \rbrace \textrm{ on } \lbrace \tau = n \rbrace, n=1,...,N.
\end{displaymath}
Define the moments with the highest number of the observed value, i.e.
\begin{displaymath}
\tau_{1}=1, \quad \tau_{k}= \inf \lbrace n: \tau_{k-1} \leq n \leq N, X_{n} = \max \lbrace X_{1},...,X_{n} \rbrace \rbrace \textrm{ for } k=1,...,N.
\end{displaymath} 
We observe that the sequence $\tau_1, \tau_2,... \in \mathcal{T}_{0}$.
Now let us consider the following chain
\begin{displaymath}
Y_{k}= (\tau_k, X_{\tau_k}) \textrm{ on } \lbrace \tau_k < N+1 \rbrace, Y_{k}=(\tau_{N+1},\xi), 
\end{displaymath}
where $\xi$ is a special absorbing state. It is easy to see that $\lbrace Y_{k} \rbrace_{k=1}^{N+1}$ is a Markov chain with transition probabilities
\begin{equation}
\begin{split}
p((n,x),(m,B)) &= P(\tau_{k+1}=m, X_{m} \in B | \tau_{k}=n, X_n=x)\\
&= x^{m-n-1}\int_{B}dy, 
\end{split}
\end{equation} 
$m>n$ and $0$ otherwise, with $B \subseteq (x,1]$.

\subsection{\label{sec2.1}Full information duration problem (FIDP) - the classical version.}
In classical FIDP we stop at the relatively best object and hold it as long as it is the relatively best one. Let $w(n,x)$ denote the expected payoff given that the $n$-th object is the relatively best object of value $X_n=x$ and we select it. It is easy to see that
\begin{equation}\label{wnx}
w(n,x) = \sum_{m=n+1}^{N+1}p((n,x),(m,B))(m-n).
\end{equation}
In our problem we assume that $B=(x,1]$ and $B=\xi$ for $m=N+1$. Therefore \eqref{wnx} has the form
\begin{equation}
w(n,x) = \sum_{m=n+1}^{N+1} x^{m-n-1} = \sum_{m=0}^{N-n}x^{m} = \dfrac{1-x^{N-n+1}}{1-x}.
\end{equation}
Let us denote $s:=N-n+1$, and $\tilde{w}(x,s):= \frac{1-x^{s}}{1-x}$. This notation we can understand as a stopping on $s$-th object from the end.
Using the backward induction method the optimal rule can be found. The optimal expected return when there are $s$ objects yet there is to be observed knowing that the present maximum value of the past observations is $x$ is defined by 
\begin{equation}\label{recursion}
v(x,s) = xv(x,s-1) + \int_{x}^{1} \max \lbrace \tilde{w}(y,s),v(y,s-1)\rbrace dy
\end{equation}
with the initial point $v(x,0) \equiv 0$. The following theorem gives the optimal stopping rule.
\begin{theorem}(Ferguson, Hardwick, Tamaki~\cite{ferhartam92:own})
In the FIDP it is optimal to select the relatively best object of value $X_s=x$ at $s$ stages from the end if $x\geq x_{s}$, where $x_{1}=0$ and for $s>1, x_s$ is the unique root of the equation
\begin{equation}\label{eqn1}
\sum_{i=1}^{s}x^{i-1} = \sum_{i=1}^{s-1}x^{i-1} \sum_{j=1}^{s-i}\frac{1-x^{j}}{j}.
\end{equation}
\end{theorem}
Since this problem is monotone then one-step-look-ahead rule is optimal. We do not write here a proof. This problem is related to the full information random horizon best choice problem of Porosi\'nski~\cite{por87:random}. The optimal rules in this two problems are the same. In \cite{por87:random} the author gives an asymptotic approximation for the $x_n$ as $n\rightarrow \infty$. Substituting $x_n :=1-\frac{z_n}{n}$ to \eqref{eqn1} and keeping this equation stay zero $z_n$ must converge to some constant $z$, where $z$ satisfies the equation 
\begin{equation}
\int_{0}^{z} e^{t}(1-\int_{0}^{t}\frac{1}{u}(1-e^{-u})du)dt=0.
\end{equation}
Using numerical methods it can be found that $z \approx 2.1198$. Therefore $x_{n} \approx 1- \frac{2.1198}{n}$. Determining the value of the problem as
$$V(n) = v(0,n) $$
in \cite{MazTam03:Explicite} there is an approximated win probability given by
$$\dfrac{V(n)}{n} \rightarrow C $$
as $n\rightarrow \infty$, where $C$ is constant given by
\begin{displaymath}
C = \int_{0}^{1} e^{-\frac{z}{u}}(\int_{0}^{u}(\dfrac{e^{\frac{zt}{u}}-1}{t}+\dfrac{e^{\frac{zt}{u}}}{1-t})dt-1)du
\end{displaymath}
using the numerical methods $C \approx 0.435171$.
The payoff function in this problem is almost the same as the payoff function in the full information best choice problem with the random, uniformly discrete distributed on $\lbrace1,...,n\rbrace$ horizon. The relation is that the payoff in FIDP is $n$ times greater than in corresponding best choice problem (see: Gnedin \cite{Gne05:Objectives}). 

In the problem where recall is possible to find the optimal strategy we use 1-Step-Look-Ahead method. If we stop at $s$ stages from the end with an object of value $x$, then our expected payoff is $\tilde{w}(x,s)$. Continuing one step we expect to receive
\begin{displaymath}
x\tilde{w}(x,s-1) + \int_{x}^{1}\tilde{w}(y,s-1)dy.
\end{displaymath}
After the simple calculation we get the optimal rule. The result is in the following
\begin{theorem}(Ferguson, Hardwick, Tamaki~\cite{ferhartam92:own})
In the FIDP with recall it is optimal to stop at $s$ stages from the end if the largest of the values observed is at least $x_s$, where $x_s \in [0,1]$ is the root of the following equation
\begin{equation}\label{eqn2}
\sum_{j=1}^{s-1}\dfrac{1-x^{j}}{j}=1.
\end{equation}
\end{theorem}

This is also the solution of the full information random horizon best choice problem with recall. Substituting $x :=1-\frac{z}{n}$ to \eqref{eqn1} and taking the limit $n\rightarrow \infty$ we need to find the solution of the following integral
\begin{equation}
\int_{0}^{1}\dfrac{1-e^{-zy}}{y}dy = 1.
\end{equation}
Numerical methods give us $x_n\approx 1-\frac{1.345}{n}$.

\subsection{\label{sec2.2}Best choice duration problem (BCDP).}
In BCDP we win the duration of owning the relatively best object only if  it is the best object overall.

Once again we introduce Markovian approach as in Section \ref{sec1}. In looking for the best object we need to stop at the last stage before the process terminates. Let $w(n,x)$ denote the expected payoff given that the $n$-th object is the relatively best object of value $X_n=x$, we select it and it is the maximum value overall. It can be observed that
\begin{equation}\label{wnx2}
w(n,x) = \sum_{m=n+1}^{N+1}p((n,x),(N+1,\xi))f(m,B).
\end{equation}
where $B=(x,1]$ and $B=\xi$ for $m=N+1$. Therefore \eqref{wnx2} has the form
\begin{equation}
w(n,x) = \sum_{m=n+1}^{N+1} x^{N-n} = x^{N-n} \sum_{m=n+1}^{N+1} = x^{N-n}(N-n+1).
\end{equation}
Let us denote $s:=N-n+1$, and $\tilde{w}(x,s):= sx^{s-1}$. Once again we can think about this notation as a stopping at $s$th candidate from the end. In this problem recursion \eqref{recursion} can be used. As a result we get the following
\begin{theorem}
In the full information BCDP it is optimal to select the relatively best object of the value $X_s=x$ at $s$ stages from the end if $x\geq x_s$, where $x_1=0$ and $x_s, s>1$ is the unique root of the equation
\begin{equation}
\sum_{j=1}^{s-1}x^{j-1} = sx^{s-1}.
\end{equation}
\end{theorem}
The value of the problem was found by Tamaki in \cite{Tam10:MonotoneThres}. It can be shown that the optimal limiting payoff is given by
\begin{equation}
v^{*}=\int_{0}^{1} \dfrac{1}{x} \left[\int_{0}^{x} e^{-\frac{c^{*}x}{1-y}}dy\right[dx - 2\int_{0}^{1}ye^{-\frac{c^{*}}{y}}dy \approx 0.31096,
\end{equation}
where $c^{*} \approx 1.2564$ is a unique solution of the equation
\begin{displaymath}
e^c=1+2c.
\end{displaymath}
With recall the optimal limiting payoff is given by
\begin{equation}
u^{*}= \dfrac{1-\log (2)}{2}+(\log (2))^2 I(\log(2)) \approx 0.33536,
\end{equation}
where $I(\cdot)$ is given by
\begin{displaymath}
I(c)=\int_{1}^{\infty}\dfrac{e^{-ct}}{t}dt.
\end{displaymath}

\subsection{\label{sec2.3}Random Horizon Full information duration problem (RHFIDP).} Let us assume that the number of actually available objects $N$ is random. This additional uncertainty in the FIDP  will cause the necessity of a clarification of the very concept of the duration. The random horizon duration problem (RHDP) for the no-information case was investigated by Tamaki~\cite{Tam13:NIDPRH}. He extended the research on RHDP for the full-information case in \cite{Tam13:RHDPFI}. Let $N$ be assumed a bounded random variable, independent of the sequence $X_1, X_2,\ldots, X_n$, and having a prior distribution $\vec{p} = (p_1, p_2,\ldots, p_n)$, where $p_k = \bP\{N = k\}$ are such that $\sum_{k=1}^n p_k=1$ and $p_n > 0$ for a known upper bound $n$. In \cite{Tam13:RHDPFI} the RHDP was distinguished into two models according to whether the final stage of the planning horizon is $N$ or $n$. This distinction is related to the last relative maximum. That is, if the chosen object is the last relative maximum prior to $N$, we hold it until the stage $N$ in the first model, whereas until the stage $n$ in the second model. The classical, finite deterministic horizon FIDP occurs as a special case of the RHDP if $N$ degenerates to $n$. In this case there is no difference between these two approaches. The performance measure is the proportional duration of holding the candidate (the relative maximum of the observed sequence).

In this paper the extension to the unbounded random horizon FIDP will be treated. For further consideration let us formulate first the main elements of the Tamaki's model with planning horizon $N$ investigated in~\cite{Tam13:RHDPFI}. For a given prior $\vec{p}$ let $\pi_k=\bP\{N\geq k\}$. Let $(k,x)$ be the state, where it is just observed, $k$-th observation to be a candidate having the value $x$, i.e $X_k=\max\{X_1,\ldots,X_k\}=x$, $1\leq k\leq n$, $0<x<1$. Denote by $s_k(x)$ the payoff earned by stopping with the current candidate in state $(k,x)$. It is the expected value of the proportional duration $D_k(x)$ when the candidate is the state $(k,x)$ (see Lemma 2.1 in~\cite{Tam13:RHDPFI}). Let $I_j(x)=\chi_{\{\omega: X_j(\omega)<x\}}(\omega)$ for $1\leq j\leq n$, $0<x<1$. The proportional duration (see Appendix A of \cite{Tam13:RHDPFI}) is 
\begin{equation}\label{D(k,x)}
D_k(x)=\frac{1}{n}\left[1+\sum_{i=k+1}^N\prod_{j=k+1}^i I_j(x)\right]
\end{equation}     
and the payoff has the form:
\begin{equation}\label{s(k,x)1}
s_k(x)=\frac{1}{n}\frac{1}{\pi_k x^k}\sum_{i=k}^n\pi_i x^i.
\end{equation}     

Also let $c_k(x)$ denote the payoff earned by continuing observations in an optimal manner. Then $v_k(x)=\max \lbrace s_k(x),c_k(x)\rbrace$ is the optimal payoff provided that we start from the state $(k,x)$. If the decision maker decides to proceed to the next stage after leaving the state $(k,x)$, the $(k+1)$-th object can be observed only with probability $\bP\{N\geq k+1|N\geq k\}=\frac{\pi_{k+1}}{\pi_k}$ and it gives the state $(k+1,y)$ if this new observation is a candidate with the value $y$, which is bigger than $x$, while it pretends to leave the state $(k+1,x)$ if it is not a candidate. It leads to the recursive equation:
\begin{equation}\label{c(k,x)}
c_k(x)=\frac{\pi_{k+1}}{\pi_k}\left[xc_{k+1}(x)+\int_x^1v_{k+1}(y)dy\right],\text{ $1\leq k<n$,}
\end{equation}
with the boundary condition $c_n(x)\equiv 0$. The repeated use of \eqref{c(k,x)} yields 
\begin{equation}\label{c(k,x)2}
c_k(x)=\sum_{i=k+1}^n\frac{\pi_{i}}{\pi_k}x^{i-k-1}\int_x^1s_i(y)dy.
\end{equation}
Since, for the given $k$, $s_k(x)$ is an increasing function of $x$ by \eqref{s(k,x)1} while $c_k(x)$ is a non-increasing function in $x$ from its definition, then there exists a sequence of thresholds $\vec{a}^\star=(a_1^\star,\ldots,a_n^\star)$ such that, when in the state $(k,x)$, the optimal rule stops with the current candidate iff $x\geq a_k^\star$ i.e. $\tau_n(\vec(a)^\star)=\min\{k: X_k=\max\{X_1,\ldots,X_k\}\geq a_k^\star\}\wedge n$. 

\subsection{\label{sec2.4} The unbounded horizon full-information duration problem.}
Let us consider the duration problem when the horizon $N$ is unbounded. The equations \eqref{D(k,x)} and \eqref{s(k,x)1} depend on the tail probabilities of the horizon distribution. They allow to formulate the equation for the value function of the problem in the form \eqref{c(k,x)} with an additional requirement that the solution $v_k(x)$ should be a continuous function of $x$ (a smooth fit condition). The methods of solving such an equation are well known in MDP theory. The special case of the horizon with the geometric distribution is worth be solved due to many reasons. Let us observe that $s_k(x)$ is independent of $k$ because the duration of the candidate at the state $(k,x)$ is $s$ by two exclusive ways:
\begin{enumerate}
\item $N\geq k+s$ and $X_{k+1}<x$, $\ldots$, $X_{k+s}>x$;
\item $N=k+s-1$ and  $X_{k+1}<x$, $\ldots$, $X_{k+s}<x$.
\end{enumerate}
It leads to 
\begin{eqnarray}\label{s(k,x)}
s_k(x)&=&\bE[n D_k(x)]=\sum_{j=1}^\infty j[(1-x)x^{j-1}\frac{\pi_{k+j}}{\pi_{k}}+x^{j-1}\frac{\pi_{k+j-1}}{\pi_{k}}]\\
\nonumber&=&[(1-x)q+p]\sum_{j=1}^\infty j(qx)^{j-1}=[(1-x)q+p]\frac{1}{(1-qx)^2}\\
\nonumber&=&\frac{1}{1-qx}
\end{eqnarray}
for the geometric distribution. There same result is given by \eqref{s(k,x)1} when the geometric horizon is applied. 

Given that the payoff function does not depend on $k$ and the lack of memory properties of the geometric distribution also the payoff earned by continuing observations in an optimal manner $c_k(x)=c(x)$ and the optimal payoff $v_k(x)=v(x)$ provided that we start from state $(k,x)$ do not depend on $k$. We have
\begin{equation}\label{c(x)}
c(x)=q[x c(x)+\int_x^1v(y)dy],
\end{equation}
which gives $c(x)=q(1-qx)^{-1}\int_x^1v(y)dy$. The optimal payoff 
\begin{equation}
v(y)=\max\{c(x),s(x)\} =(1-qx)^{-1}\max\{1,q\int_x^1 v(y)dy\}.
\end{equation}
For $x$ enough close to $1$ we have $v(x)=s(x)$. It means that the stopping region contains $\{(k,x):x\geq x_0\}$, where $x_0$ fulfils the condition $\int_{x_0}^1v(y)dy=\frac{1}{q}$. For $x\leq x_0$ the optimal payoff fulfil eqution:
\begin{equation}\label{contFIT}
(1-qx)v(x)=a[\int_{x}^{x_0}v(y)dy+\int_{x_0}^1\frac{1}{1-qy}dy].
\end{equation}
It implies that $v(x)=\text{const}$ for $x\in(0,x_0]$. The continuity condition forces $v(x_0)=(1-qx_0)^{-1}$. If such $x_0\in (0,1)$ exists then by \eqref{contFIT} we get $\ln\frac{1-q}{1-qx_0}=-1$. The function $\varphi(t)=1+\ln\frac{p}{t}$ is well defined for $t\in (p,1)$. It is non-increasing in this domain, $\varphi(p)=1$ and $\varphi(1)<0$ when $p<\exp(-1)$. 
\begin{conclusion}
If $p\leq \exp(-1)$, then there is $x_0=\frac{1-e p}{q}\in (0,1)$ such that $\{(k,x):x\geq x_0\}$ is the optimal stopping region for RHFIDP with the geometric horizon. The expected optimal payoff is $v^\star=\frac{x_0}{1-qx_0}-\frac{1}{q}\ln\frac{1-q}{1-qx_0}$.

If $p>\exp(-1)$ then the optimal stopping region for RHFIDP is the hole state space. The decision maker should stop at the first observation obtaining the expected payoff $v^\star=-\frac{1}{q}\ln(1-q)$. 
\end{conclusion} 

\begin{remark}
The definition of the duration depends on the context. In the seminal paper by Ferguson et al.~\cite{ferhartam92:own} there are various models. In most of them the maturity of the accepted object is the moment when it stops to be the candidate. There are other cases when the maturity is related to the approach of the horizon. In the finite horizon case it is assumed that the duration is expanded by adding $1$. When the horizon is random, as in this section, the same understanding of the maturity is applied. However, in various applications such a definition of the maturity should be corrected. If we buy obligation then it has an additional value to the prescribed moment. If this maturity moment is random we can consider the case when the random horizon is observed immediately or we learn about the maturity by symptoms like an absence of new observations. 

Assuming the maturity as an immediate close when the decision maker reaches the last observation the expected payoff defined by \eqref{s(k,x)} is changed to 
\begin{eqnarray}\label{s(k,x)2}
\stilde_k(x)&=&\bE[n \Dtilde_k(x)]=\sum_{j=1}^\infty j[(1-x)x^{j-1}\frac{\pi_{k+j}}{\pi_{k}}+x^{j}\frac{\pi_{k+j}}{\pi_{k}}]\\
\nonumber&=&[(1-x)q+pqx]\sum_{j=1}^\infty j(qx)^{j-1}=[(1-x)q+pqx]\frac{1}{(1-qx)^2}\\
\nonumber&=&\frac{q}{1-qx}.
\end{eqnarray}

The optimal strategy does not change with respect to the previous model but the expected optimal payoff does.
\end{remark}

\subsection{\label{sec2.4a}Duration of owning relatively best or second best object.}
%\input{sec2.4a.tex}
%%sec2.4a
This problem was firstly considered by Kurushima and Ano. The objective is to maximize the time period of owning the relatively best and the relatively second-best object. Here we consider the class of the stopping rule restricted only to stop at the relatively best object. Let $U_{n}(x)$ denote the expected duration of the relatively best object whose rank remains within the two when the time to go is $n$ and the decision maker accepts the relatively best applicant whose value $x$ is the maximum value among that of the applicants arrived so far, that is $X _{n}=x$. $U_{n}(x)$ is given by
\begin{equation}
U_{n}(x)=2\sum_{k=1}^{n-1}x^{k-1}-nx^{n-1}.
\end{equation}
Let $G_{n}(x)$ define 
\begin{equation}
G_{n}(x)=U_{n}(x)- \sum_{k=1}^{n-1}x^{k-1}\int_{x}^{1}U_{n-k}(y)dy.
\end{equation}
The 1-SLA calls for a stop in the region $B$ which is described as
$B=\lbrace(n, x):G_{n}(x)\geq 0\rbrace$, where $(n, x)$ is represented by the state when the time to go is $n$ and the present applicant is the relatively best one whose value $x$ is the maximum value among the applicants arrived so far.
To show that the 1-SLA stopping rule is optimal, it is sufficient to show the next two statements: 
\begin{enumerate}
\item $G_{n}(x)\geq 0\Rightarrow G_{n-k}(x)\geq 0, k=1,2, \ldots$
\item $G_{n}(x)\geq 0\Rightarrow G_{n}(y)\geq 0, y\geq x$.
\end{enumerate}
The first statement is presented in the mentioned paper. The second remains as an open problem. This problem is concluded by
\begin{conjecture} (see: \cite{KurAno10:FIDP})
For the full-information case of the duration problem where the objective is to maximize the duration of owning the relatively best or the second-best object, we assume that the class of the stopping rule is restricted to that of stopping only at the relatively best object. Then, the optimal stopping rule is to accept the first applicant who has the maximum $X_{n}=x\geq s_{n}$  among the observed objects so far when the remaining time is $n$, where $s_{1}=1$ and $s_{n}, n\geq 2$ is the unique root of the equation
$$
3\sum_{k=1}^{n}x^{k-1}-2nx^{n-1}-2\sum_{k=1}^{n-1}x^{k-1}\sum_{j=1}^{n-k-1}\frac{1}{j}+2\sum_{k=1}^{n-1}x^{k}\sum_{j=1}^{k}\frac{1}{j}=0.$$

\end{conjecture}

\subsection{\label{sec2.5}Duration of owning relatively best or second best object for unbounded horizon.}
\subsubsection{Applying geometrical horizon in Kurushima and Ano problem.}
Our aim is to maximize the duration of owning the relatively best or the second best object, where the class of stopping times is restricted to the relatively best objects. We observe $N$ random variables from the known distribution. We consider a special case where $N$ is a random variable geometrically distributed, i.e.
\begin{equation}\label{geom}
P(N=k) = p_k = pq^{k-1}, \quad 0 < p < 1; q=1-p.
\end{equation}
Let $w(n,x)$ denote the payoff for stopping at the $n$-th object whose value is $x$. 
\begin{equation}\label{equation2.5.1}
w(n,x) = \bE[T_n-n|X_n=x]= \frac{2q}{1-qx}-\frac{q(1-q)}{(1-qx)^2}.
\end{equation}
If we continue observations we expect to get a reward $Tw(n,x)$:
\begin{equation}\label{equation2.5.2}
\begin{split}
Tw(n,x) &=\sum_{m=n+1}^{\infty}\int_{x}^{1}w(n,v)\dfrac{\pi_m}{\pi_n}x^{m-n-1}dv \\
&= \dfrac{q}{1-qx}\Big( 2\log\big( \frac{qx-1}{q-1} \big) + \frac{1-q}{1-qx}-1 \Big).
\end{split}\end{equation} 
Denote $G(n,x):=w(n,x)-Tw(n,x)$. The 1-SLA rule is described by the following set:
$$ B=\lbrace (n,x): G(n,x) \geq 0 \rbrace. $$
To prove the optimality of the 1-SLA rule it is necessary to show two things: (i)  $G(n,x)\geq 0 \Rightarrow G(n+k,x)\geq 0, k=1,2,...$ and (ii) $G(n,x)\geq 0 \Rightarrow G(n,y), y\geq x$. (i) is obvious, because payoffs do not depend on $n$. Consider the function $g(x):=G(n,x)$. We will show that if $g(x)\geq x$ for some $x$ then it will be greater than $0$ for $y>x$. It is an equivalent to the statement if for a given $x$ inequality
$$3-2\dfrac{1-q}{1-qx} - \log \big( \frac{1-qx}{1-q} \big)^{2} \geq 0, $$
holds then the same is for every $y>x$. Calculating a derivative of LHS we get
%\begin{displaymath}
%\begin{split}
\[ -2q(1-q)\dfrac{1}{(1-qx)^2}+2q(1-qx)\dfrac{1}{(1-qx)^2} = \dfrac{2q^2(1-x)}{(1-qx)^2}.
\]
%\end{split}
%\end{displaymath}
It is non-negative for every $x\leq 1$ and in point $x=1$ it has the value $1$. If the inequality on LHS is true for the given $x$ then it is true for $x<y<1$. We conclude that the 1-SLA rule is optimal. It is a threshold strategy with threshold
\begin{equation}\label{prog}
\dfrac{\mu^{\star}p-1}{p-1},
\end{equation}
where $\mu^{\star}$ is a solution of the equality
$$ \mu^2e^{\frac{2}{\mu}}=e^3, \mu>1.$$
Numerically $\mu^{\star} \approx 3.3145.$ Note that if $p\geq \frac{1}{\mu^{\star}} \approx 0.3017046 $ the threshold is $0$ so the 1-SLA rule calls for a stop at the very first observation.
%%%
\subsubsection{Stopping on the best and the second best object.} 
Assume that at the moment $n$ we are interested in choosing the relatively best or second best object. We choose the object and hold it as long as it is the relatively best or second best object. 
We observe sequentially $X_1,...,X_N$ i.i.d. $N$ is as in \eqref{geom}. Let $\zeta_n$ denote the largest observation of sequence $X_1,...,X_n$ and similarly $\eta_n$ denote the second largest value of the sequence ($\zeta_{1}=X_{1}, \zeta_{n}=X_{n:n}$ for $ n>1,  \eta_{1}=0, \eta_{n}= X_{n-1:n} $ for $ n>1$). Let $T_n$ be a random variable that denotes the moment after the time $n$ when a better observation occurs. Our aim is to find such a stopping time $\tau^{*} \in \mathcal{T}$ that 
$$P(\tau^{*} \leq N, T_{\tau^{*}}-\tau^{*}, X_{\tau^{*}}=\zeta_{\tau^{*}} \textrm{ or } X_{\tau^{*}}=\eta_{\tau^{*}})$$
is maximized, where $\mathcal{T}$ denotes a set of all stopping times with respect to the family $\lbrace \mathcal{F}_{n} \rbrace_{n=1}^{\infty}$. Denote $R_n$ the rank of $n$th observation, i.e. a random variable:
$$R_{n} = \begin{cases}
1,  X_n=\zeta_n\\
2,  X_n=\eta_n.
\end{cases} $$
Let
\begin{displaymath}
\begin{split}
\tau_1 & = 1 \\
\tau_{k+1} & = \inf \lbrace n: \tau_{k}<n\leq N, X_{n} \geq \eta_{\tau_{k}}, k\in\mathbb{N} \rbrace.
\end{split}
\end{displaymath}
Let us define a sequence
\begin{equation}\label{chain}
\begin{split}
Y_{k}&=(\tau_k,\zeta_{\tau_k}, \eta_{\tau_k}, R_{\tau_k}), \quad \textrm{ if } \tau_k<\infty \\
Y_{k}&= \delta, \quad \textrm{ if } \tau_k = \infty,
\end{split}
\end{equation}
where $\delta$ is the special absorbing state. It is easy to verify that \eqref{chain} is a Markov chain with respect to $\mathcal{F}_{\tau_k}= \sigma (X_{1},...,X_{\tau_k}, \mathbb{I}_{(N\geq 0)}, ..., \mathbb{I}_{(N\geq \tau_k-1)}), k=1,2... .$ The transition probabilities are:
\begin{equation}\label{transition1}
p((n,x,y,i),(m,[0,u],x,1))= \begin{cases}
\dfrac{\pi_m}{\pi_n}y^{m-n-1}(u-x), u>x \\
0, \textrm{ otherwise},
\end{cases}
\end{equation}
\begin{equation}\label{transition2}
p((n,x,y,i),(m,x,[0,u],2))= \begin{cases}
\dfrac{\pi_m}{\pi_n}y^{m-n-1}(u-y), y\leq u \leq x \\
%\dfrac{\pi_m}{\pi_n}y^{m-n-1}(x-y), u \geq x\\
0,  \textrm{ otherwise},
\end{cases}
\end{equation}
for $m>n, i=1,2$, where $\pi_{k}=\sum_{j=k}^{\infty}p_j = q^{k-1}$. We derive the payoff function. If we stop at the state $(n,x,y,1)$ then the duration of the candidate is $k$ if:
\begin{enumerate}
\item if there are two candidates of rank $1$: one at the time $n+1$ and the second at some point between $n+1$ and $n+k-1$ inclusive and the time horizon is longer than $n+k$  
\item if there is one candidate of rank $1$ between $n+1$ and $n+k-1$ inclusive and the time horizon is $n+k$ 
\item if there is no more better candidates and time horizon is $n+k$ 
\end{enumerate}
If we stop at state $(n,x,y,2)$ then the duration of candidate is $k$ if:
\begin{enumerate}
\item if there is one candidate of rank $1$ or $2$ at the time $n+k$ and the time horizon is longer than $n+k$ 
\item if there are no more better candidates and the time horizon is $n+k$
\end{enumerate} The payoff function is given by:
\begin{equation}\label{payoff}
\begin{split}
W((n,x,y,1)) &= \bE[T_n-n|X_n=x,R_n=1,\zeta_n=x,\eta_n=y] \\
& = \frac{2q}{1-qx}-\frac{q(1-q)}{(1-qx)^2}\\
W((n,x,y,2)) &= \bE[T_n-n|X_n=y,R_n=2,\zeta_n=x,\eta_n=y] \\
& = \frac{q}{1-qy}
\end{split}
\end{equation}
Note that both functions do not depend on $n$. This will not be mentioned later. Let $T$ be an operator $Tf(x) = \int_{X}f(z)dP_{x}(z)$ for the bounded function $f:X\rightarrow\mathbb{R}$. Then
\begin{equation}
\begin{split}
TW(x,y,i) & =\dfrac{q}{1-qy}\Big( \log\big(\frac{(1-qy)(1-qx)}{(1-q)^2} \big) + \frac{1-q}{1-qx}-1 \Big)\quad i=1,2.
\end{split}
\end{equation}
We transform our equations in the following way:
\begin{displaymath}
s:=\frac{1-q}{1-qx}, \qquad t:=\frac{1-q}{1-qy}, \qquad \alpha :=\frac{q}{1-q}, 
\end{displaymath}
and $\frac{1}{1+\alpha} \leq t \leq s \leq 1$, where $\alpha \in [0,\infty) $. Our payoff function is now given by:
\begin{equation}
\begin{split}
W(s,t,1) &= \alpha s(2-s)\\
W(s,t,2) &= \alpha t\\
TW(s,t,i)&= \alpha t\left( -\log(st)+s-1 \right), \quad i=1,2.
\end{split}
\end{equation}
Let us denote $F(s,t):= W(s,t,1)-TW(s,t,1)$ and $G(s,t):= W(s,t,2)-TW(s,t,2)$.
To find the stopping set we need to find the set that satisfies the optimality equation. Let $B_1 = \lbrace (s,t): F(s,t) \geq 0 \rbrace, B_2 = \lbrace (s,t): G(s,t) \geq 0 \rbrace$. \\

Note that the sign of the functions $F$ and $G$ does not depend on the parameter $\alpha$. These sets are the same as in \cite{porsza90:selection} (page 685). Using the same methodology we obtain the optimal strategy: it allows us to stop only in such a moment, when observed candidate is the largest so far and exceeds the value $x^{\star}$. This value is independent of the value of the second largest object. The problem is now transformed into the problem of stopping only at the largest observations, as it was in Kurushima and Ano problem with the geometrical horizon. Even if we observe the value of the second largest object it does not affect the optimal strategy. We can reduce the problem to observing only the relatively best objects.

%*******************************************************************
%THIRD SECTION - DO NOT FORGET THE SETCOUNTER{EQUATION}{0}
%*******************************************************************

%\vspace*{0.5cm}
%\setcounter{equation}{0}
%\subsection{Discounted full information duration problem (DFIDP)}\label{sec3}

%*******************************************************************
%FOURTH SECTION - DO NOT FORGET THE SETCOUNTER{EQUATION}{0}
%*******************************************************************

%*******************************************************************
%FIFTH SECTION - DO NOT FORGET THE SETCOUNTER{EQUATION}{0}
%*******************************************************************

%*******************************************************************
%CONCLUSIONS - DO NOT FORGET THE SETCOUNTER{EQUATION}{0}
%*******************************************************************

%\vspace{-1ex}
\section{Conclusion and acknowledgements}\label{koniec}
We have presented different results regarding the duration problem. This is the very first time, when all results from different papers are collected in one place. The previous results were presented in a new form and thanks to that it is easy to see connection between two problems. A very intriguing issue seems to be the problem with random geometrically distributed horizon. The geometric distribution has memoryless property and this is the reason why the result is quite simple. \\

We are grateful to an anonymous referee for the several comments and improvements of this paper.
%The authors would like to thank the referee for advices.

%\section{Acknowledgements}
%We are grateful to an anonymous referee for several comments and improvements of this paper.

\nocite{MazTam07:Trajectories}

%*******************************************************************
%BIBLIOGRAPHY
%*******************************************************************
%%\bibliography{compet,compet3,prace2}
%%\bibliographystyle{abbrv}

\end{document}